\documentclass[smallextended]{svjour3}

\smartqed
\usepackage{hyperref}
\usepackage{amsmath}
\usepackage{amsmath,amssymb}
\usepackage{xfrac}
\usepackage{microtype}
\usepackage{setspace}
\usepackage[english]{babel}
\usepackage{cite}
\usepackage{graphicx}
\usepackage{caption}
\usepackage{subcaption}
\usepackage{xcolor}
\usepackage{bm}
\DeclareGraphicsExtensions{.pdf}
\journalname{Anonymous}
\begin{document}

\title{Beyond the stability}
\titlerunning{Beyond the stability}
\author{Majid Akbarian, University of Bojnord, Iran}

\maketitle
\

\section{Discussion}
\label{Sec. 1}

Consider
\begin{equation}\label{equ.2}
	\begin{cases}
		\dot{x} = f(t,x), \\
		x(t_0) = x_0,
	\end{cases}
\end{equation}
where $f:[0,\infty ]\times D\to\mathbb{R}^{n}$ is locally Lipschitz in $x$ on {${[}{0}{,}{\infty }{]\times }{D}$}, and piecewise continuous in $t$. Also, $D \subset \mathbb{R}^n$ is an open set containing $x=0$. The origin is the equilibrium point for system (\ref{equ.2}) and $t_0$ is the initial time.\\

We recall from \cite{7} that the origin of (\ref{equ.2}) is uniformly stable when\\ 
\begin{equation}\label{equ.3}
	\begin{aligned}
		&\text{for all }\varepsilon > 0, \quad \text{there exists }\delta =\delta (\varepsilon )>0, \quad s.t.\ \|x_0\|<\delta\Rightarrow \|x(t)\|<\varepsilon \\
		&\text{for all } t\ge t_0\ge 0.
	\end{aligned}
\end{equation}
In fact, if the equilibrium point of (\ref{equ.2}) is uniformly stable, $\delta $ does not depend on $t_0$ explicitly.\\

Additionally, the origin of equation (\ref{equ.2}) exhibits the following stability properties:
	
	\begin{itemize}
		\item \textbf{Uniformly asymptotically stable} if it is uniformly stable and there exists a positive constant $c$, independent of $t_0$, such that for all initial conditions $\lVert x(t_0) \rVert < c$, the solution $x(t)$ approaches zero as $t \to \infty$, with uniformity over $t_0$. Specifically, for every $\eta > 0$, there exists a time $T = T(\eta) > 0$ satisfying
		\begin{equation}\label{equ.4}
			\lVert x(t) \rVert < \eta, \quad \forall t \geq t_0 + T(\eta), \quad \forall \lVert x(t_0) \rVert < c 
		\end{equation}
		\item \textbf{Globally uniformly asymptotically stable} if it is uniformly stable, and $\delta(\epsilon)$ can be selected such that $\lim_{\epsilon \to \infty} \delta(\epsilon) = \infty$. Moreover, for any pair of positive numbers $\eta$ and $c$, there exists a time $T = T(\eta, c) > 0$ such that
		\begin{equation}\label{equ.5}
			\lVert x(t) \rVert < \eta, \quad \forall t \geq t_0 + T(\eta, c), \quad \forall \lVert x(t_0) \rVert < c 
		\end{equation}
	\end{itemize}

	\textbf{Theorem 1} Let $x = 0$ be an equilibrium point for equation (\ref{equ.2}), and let $D \subset \mathbb{R}^n$ be a domain containing $x = 0$. Suppose there exists a pair of functions $(V, W^*)$ where $V : [0, \infty) \times D \to \mathbb{R}$ is a continuously differentiable function satisfying
	\begin{equation}\label{equ.6}
		V_1(x) \leq V(t, x) \leq V_2(x),
	\end{equation}
	and $W^* : [0, \infty) \times D \to \mathbb{R}$ is a continuous function with $W^*(t, 0) = 0$, such that
	\begin{equation}\label{equ.7}
		\int_0^\infty \max \{ W^*(\tau, x(\tau)), 0 \} \, d\tau \leq M(x_0),
	\end{equation}
	where $M(x_0)$ is a bounded and continuous function of $x_0$. Additionally, the following holds:
	\begin{equation}\label{equ.8}
		\frac{\partial V}{\partial t} + \frac{\partial V}{\partial x} f(t, x) - \max \{ W^*(t, x), 0 \}\leq 0,
	\end{equation}
	for all $t \geq 0$ and $x \in D$, with $V_1(x)$ and $V_2(x)$ being continuous positive definite functions on $D$. Under these conditions, $x = 0$ is uniformly stable.\\
	\\
\textbf{Theorem 2} Assume that the conditions of Theorem 1 hold, with inequality (\ref{equ.8}) strengthened to
\begin{equation}\label{equ.9}
		\frac{\partial V}{\partial t} + \frac{\partial V}{\partial x} f(t, x) - \max \{ W^*(t, x), 0 \} \leq -V_3(x),
\end{equation}
for all $t \geq 0$ and $x \in D$, where $V_3(x)$ is a continuous positive definite function on $D$. Then, the equilibrium point $x = 0$ is uniformly asymptotically stable. Furthermore, if $D = \mathbb{R}^n$ and $V_1(x)$ is radially unbounded, then $x = 0$ is globally uniformly asymptotically stable.\\

%%%%%%%%%%%%%%%%%%%%%%%%%%%%%%%%%%%%%%%%5
\textbf{Definition 1}\cite{8} We say that $\dot{W}(x, t)$ is definitely not equal to zero in the ensemble $E(V^* = 0)$ if for any numbers $\alpha$ and $A$, $(0 < \alpha < A < H)$ numbers $r_1(\alpha, A)$, $\xi(\alpha, A)$ $(0 < r_1 < d)$, $\xi > 0$, can be found such that $|\dot{W}(x, t)| > \xi$ for $\alpha < \lVert x \rVert < A$, $\rho(x, E(V^* = 0)) < r_1$, $t \geq 0$ where
\[
\rho(x, x^0) = \sqrt{(x_1 - x_1^0)^2 + \dots + (x_n - x_n^0)^2}
\]
\[
\rho(x, M) = \inf \{ \rho(x, x^0), x^0 \in M \}
\]

is called distance of the point \( x \) to the point \( x^0 \) in \( \mathbb{R}^n \) (corresponding to the elements \( M \subset \mathbb{R}^n \)).\\

\textbf{Theorem 3}\cite{8} Let there be functions $V(x, t)$, $\dot{W}(x, t)$, having in $D$ the following properties:

\begin{enumerate}
	\item The function $V(x, t)$ is positive definite and admits an infinitely small upper limit.
	\item The derivative $\dot{V}(x, t) \leq V^*(x) \leq 0$.
	\item The function $W(x, t)$ is bounded.
	\item $\dot{W}(x, t)$ is definitely not equal to zero in the ensemble $E(V^* = 0)$.
\end{enumerate}

Then the origin of system (\ref{equ.2}) is asymptotically stable with respect to $x_0$, $t_0$.\\
%%%%%%%%%%%%%%%%%%%555

\textbf{Corollary 1} Consider the non-autonomous system described by (\ref{equ.2}). For this systems, if the assumptions of Theorem 3 is satisfied, then, there exists a pair of $(V,W^*)$ that satisfy the conditions of Theorem 2. \\

\textbf{Proof.} According to Assumptions 1 and 2 in Theorem 3, the origin of system~(\ref{equ.2}) is uniformly stable. Consequently, by Lemma 4.5 in~\cite{7}, there exist a class-$\mathcal{K}$ function $\alpha_1$ and a positive constant $c$ such that
\begin{equation}\label{equ.10}
	\|x(t)\| \leq \alpha_1(\|x_0\|), \quad \text{for all } t \geq 0 \text{ and } \|x_0\| < c.
\end{equation}
Since $\dot{W}(t,x)$ is continuous and bounded, Lemma 4.3 in~\cite{7} implies that for any ball $B_r \subset D$ with $r > 0$, there exists a class-$\mathcal{K}$ function $\alpha_2$ such that
\begin{equation}\label{equ.11}
	W(t,x(t)) \leq \alpha_2(\|x(t)\|) \leq \alpha_2(\alpha_1(\|x_0\|)) =: \alpha_3(\|x_0\|).
\end{equation}

Now, let the function \( V(t,x) \) in Theorem 2 be the same as that given in Theorem 3, and define
\begin{equation}\label{equ.12}
	W^*(t,x) =
	\begin{cases}
		0, & V^*(x) < 0,\\
		|\dot{W}(t,x)|, & V^*(x) = 0.
	\end{cases}
\end{equation}
Suppose \( x(t) \in E := \{x \in D \mid V^*(x) = 0\} \). We claim that
\[
\int_0^\infty \max\{W^*(\tau, x(\tau)), 0\} \, d\tau = \int_{t_0}^\infty |\dot{W}(x(\tau), \tau)| \, d\tau
\]
is bounded by a continuous function of \( x_0 \).

Define the set
\begin{equation}\label{equ.13}
	U := \left\{ x \in D \,\middle|\, \alpha < \|x\| < A,\; \rho(x, E) < r_1 \right\},
\end{equation}
where \( \rho(x, E) \) denotes the distance from \( x \) to the set \( E \). 

Following a similar argument to that in the proof of Theorem 1.1 in~\cite{8}, any trajectory of system~(\ref{equ.2}) starting in \( E \) cannot remain in \( U \) for an interval of time longer than \( \frac{2L}{\xi} \), where \( |\dot{W}| > \xi \) and \( |W| \leq L \) in \( U \).

Let \( \{[t_k, t_k + \Delta t_k]\}_{k=1}^N \) denote the collection of time intervals where \( x(t) \in U \). Then by Theorem 1.1, each interval satisfies \( \Delta t_k \leq \frac{2L}{\xi} \), and by~(\ref{equ.11}),
\begin{equation}\label{equ.14}
	\int_{t_k}^{t_k + \Delta t_k} |\dot{W}(x(t), t)| \, dt \leq 2\alpha_3(\|x_0\|).
\end{equation}
Summing over all such intervals, we obtain:
\begin{equation}\label{equ.15}
	\int_{t_0}^\infty |\dot{W}(x(\tau), \tau)| \, d\tau \leq \sum_{k=1}^N \int_{t_k}^{t_k + \Delta t_k} |\dot{W}(x(\tau), \tau)| \, d\tau \leq N \cdot 2\alpha_3(\|x_0\|).
\end{equation}

We now claim that \( N \) must be finite. Suppose, for the sake of contradiction, that \( N \to \infty \). Using similar reasoning as in~\cite{8}, every time the trajectory enters \( U \), the value of \( V(t,x) \) decreases by at least
\[
a := \frac{\xi r_1}{2X\sqrt{n}} > 0.
\]
Thus, after \( k \) entries, we have
\[
V(x(t_k), t_k) \leq V(x_0, t_0) - ka.
\]
Letting \( k \to \infty \), this implies \( V(x(t_k), t_k) \to -\infty \), which contradicts the boundedness of \( V \). Hence, \( N \) must be finite.

Therefore, from~(\ref{equ.15}), it follows that \( \int_0^\infty \max\{W^*(\tau, x(\tau)), 0\} \, d\tau \) is bounded by a continuous function of \( \|x_0\| \).

Finally, consider the expression \( \dot{V}(t,x) - \max\{W^*(t,x), 0\} \), where
\[
\dot{V}(t,x) = \frac{\partial V}{\partial t} + \frac{\partial V}{\partial x} f(t,x).
\]
We then have:
\begin{equation}\label{equ.16}
	\dot{V}(t,x) - \max\{W^*(t,x), 0\} \leq 
	\begin{cases}
		V^*(x), & V^*(x) < 0,\\
		-|\dot{W}(t,x)|, & V^*(x) = 0,
	\end{cases}
\end{equation}
which shows that \( \dot{V}(t,x) - \max\{W^*(t,x), 0\} \) is negative definite if we define
\[
V_3(x) :=
\begin{cases}
	-V^*(x), & V^*(x) < 0,\\
	|\dot{W}(t,x)|, & V^*(x) = 0.
\end{cases}
\]
This completes the proof.
$~~~~~~~~~~~~~~~~~~~~~~~~~~~~~~~~~~~~~~~~~~~~~~~~~~~~~~~~~~~~~~~~~~~~~~~~~~~\blacksquare$\\

Here, by applying Theorem 1, we analyze the stability properties of a system considered in~\cite{9}, described as:
\begin{equation}\label{equ.17}
	\dot{x}_i = \frac{\beta(t)x_i}{1 + h(x_i)}, \quad i = 1, \dots, n,
\end{equation}
where \( \beta(t) > 0 \) with \( \int_{t_0}^\infty \beta(\tau) \, d\tau = M_1 < \infty \), and \( h(x_i) \geq 0 \). Although this system is known to be uniformly stable, a smooth Lyapunov function of the form \( V(x) \) cannot be constructed to establish its stability via classical Lyapunov methods, as shown in~\cite{9}.

However, Theorem 1 allows us to analyze the stability of system~(\ref{equ.17}) using the simple Lyapunov candidate \( V(x) = \frac{1}{2} \|x\|^2 \). Starting from~(\ref{equ.17}), we obtain:
\begin{gather}\label{equ.18}
	\begin{aligned}
		\frac{dx_i}{x_i} &= \frac{\beta(t)\, dt}{1 + h(x_i)} \\
		\Rightarrow \int_{x_i(0)}^{x_i(t)} \frac{dx_i}{x_i} &= \int_{t_0}^t \frac{\beta(\tau)}{1 + h(x_i(\tau))} \, d\tau \\
		&\leq \int_{t_0}^t \beta(\tau) \, d\tau \leq \int_{t_0}^\infty \beta(\tau) \, d\tau = M_1, \\
		\Rightarrow \ln |x_i(t)| - \ln |x_i(0)| &\leq M_1 \\
		\Rightarrow \ln \left( \frac{|x_i(t)|}{|x_i(0)|} \right) &\leq M_1 \\
		\Rightarrow |x_i(t)| &\leq |x_i(0)| e^{M_1}.
	\end{aligned}
\end{gather}
This inequality holds for all \( i = 1, \dots, n \). Hence, we have:
\begin{equation}\label{equ.19}
	\|x(t)\|^2 \leq \sum_{i=1}^n x_i^2(0) e^{2M_1} = \|x_0\|^2 e^{2M_1}.
\end{equation}

Now, using the Lyapunov function \( V(x) = \frac{1}{2} \|x\|^2 \), we compute its time derivative:
\begin{equation}\label{equ.20}
	\dot{V}(x) = \sum_{i=1}^n \frac{\beta(t) x_i^2}{1 + h(x_i)} > 0.
\end{equation}

To apply Theorem 1, define
\[
W^*(t,x) = \sum_{i=1}^n \frac{\beta(t) x_i^2}{1 + h(x_i)}.
\]
Using~(\ref{equ.19}), we estimate the integral of \( W^* \) as:
\begin{equation}\label{equ.21}
	\begin{aligned}
		\int_0^\infty \max \{ W^*(\tau, x(\tau)), 0 \} \, d\tau &\leq \int_0^\infty \sum_{i=1}^n \beta(\tau) x_i^2(\tau) \, d\tau \\
		&\leq \int_0^\infty \beta(\tau) \|x_0\|^2 e^{2M_1} \, d\tau \\
		&= \|x_0\|^2 e^{2M_1} \int_0^\infty \beta(\tau) \, d\tau = M_1 \|x_0\|^2 e^{2M_1}.
	\end{aligned}
\end{equation}

Therefore, \( W^*(t,x) \) satisfies condition~(\ref{equ.7}) in Theorem 1. Furthermore, we observe that
\[
\dot{V}(x) - \max\{W^*(t,x), 0\} = 0.
\]
This confirms, by Theorem 1, that the origin of system~(\ref{equ.17}) is uniformly stable—even though, as shown in~\cite{9} (Theorem 3), no smooth Lyapunov function of the form \( V(x) \) exists that can directly prove stability through classical methods.

\bibliographystyle {ieeetr}
\bibliography {control}

\end{document}